\documentclass{article}%
\usepackage[utf8x]{inputenc}
\usepackage{amsfonts}
\usepackage{amsmath}
\usepackage{amssymb}
\usepackage{graphicx}
\usepackage{xcolor}%
\usepackage{enumitem}

\usepackage{etoolbox}
\usepackage{hyperref}

\makeatletter
\newcounter{author}
\renewcommand*\author[1]{%
	\stepcounter{author}%
	\ifnum\c@author=1
	\gdef\@author{#1}%
	\else
	\xdef\@author{\unexpanded\expandafter{\@author\and#1}}%
	\fi
	\csgdef{author@\the\c@author}{#1}}
\newcommand*\email[1]{%
	\csgdef{email@\the\c@author}{#1}}
\newcommand*\address[1]{%
	\csgdef{address@\the\c@author}{#1}}
\AtEndDocument{%
	\xdef\author@count{\the\c@author}%
	\c@author=1
	\print@authors}
\newcommand*\print@authors{%
	\ifnum\c@author>\author@count
	\else
	\print@author{\the\c@author}%
	\advance\c@author by 1
	\expandafter\print@authors
	\fi}
\newcommand*\print@author[1]{%
	\par\medskip
	\begin{tabular}{@{}l@{}}%
		\textsc{\csuse{author@#1}}\\
		\csuse{address@#1}\\
		\textit{E-mail}:
		\href{mailto:\csuse{email@#1}}{\csuse{email@#1}}
\end{tabular}}
\makeatother

\newtheorem{theorem}{Theorem}

\newtheorem{proposition}[theorem]{Proposition}

\begin{document}

\title{Stochastic Navier-Stokes Equations and Related Models}
\author{Luigi Amedeo Bianchi}
\address{Università di Trento}
\email{luigiamedeo.bianchi@unitn.it}
\author{Franco Flandoli}
\address{Scuola Normale Superiore}
\email{franco.flandoli@sns.it}
\date{\textit{in honour of Giuseppe Da Prato}}
\maketitle

\begin{abstract}
Regularization by noise for certain classes of fluid dynamic equations, a
theme dear to Giuseppe Da Prato \cite{DaDeb2003JdMPeA}, is reviewed focusing on 3D
Navier-Stokes equations and dyadic models of turbulence.

\end{abstract}

\section{Introduction}

This is a review paper dealing with a specific question of stochastic fluid dynamics
which occupied many years of research of Giuseppe Da Prato, prepared on
the occasion of his 80th birthday. The question is whether noise may improve
the theory of well posedness of certain equations of fluid dynamics, first of
all the 3D incompressible Navier-Stokes equations.

As better remarked below, the deterministic theory of such equations has been frozen
for many years in the following duality (up to numerous side results, including very
advanced ones, which however do not change this simplified picture):
\begin{enumerate}[label=\roman*.]
	\item Weak solutions exist, globally in time, but their uniqueness is an open
	problem;
	\item More regular solutions exists uniquely locally in time, but
	their blow-up or persistence is an open problem (see the Millennium Prize
	problem described by Fefferman in~\cite{Fef2006TMPP}).
\end{enumerate} 
Why should the presence of a noise improve
such results? Specifically for the 3D Navier-Stokes equations we do not have a
precise intuition, except for the vague feeling that disorder, created
intrinsically by turbulence or imposed from outside by a noise, could
disgregate well prepared configurations which could otherwise blow up. Even if the
intuition is poor, the question is meaningful, having in mind analogous results
of \textit{regularization by noise} holding for several classes of stochastic
differential equations, ranging from classical finite dimensional cases such as~\cite{Ver1981MUS, KryRoc2005PTRF} to infinite dimensional ones, although the latter have been proved
until now only for nonlinear systems much simpler than 3D Navier-Stokes
equations, as for example~\cite{Gyo1998SPatA, DaFla2010JoFA, DaFlaPriRoc2013AP, DaFlaRocVer2016AP}. These results
prove uniqueness for certain equations with nondegenerate additive noise, in
cases where the same equations without noise miss uniqueness; and, for the
purpose of the upcoming discussion, let us mention that all of them (with the exception of~\cite{Gyo1998SPatA}) are based on suitable regularity results for the Kolmogorov
equations associated to the stochastic equations. In response to such results, the hope of proving uniqueness of weak solutions to the 3D Navier-Stokes equations
by adding a nondegenerate noise with suitable covariance rose high. Giuseppe Da Prato
made a tremendous contribution to answering this question, although for the time being the final question
is still open: together with Arnaud Debussche, in the paper~\cite{DaDeb2003JdMPeA} he constructed a smooth
solution of the infinite dimensional Kolmogorov equation associated to the
stochastic 3D Navier-Stokes equations, with a really original
and highly non trivial procedure. Existence of sufficiently smooth solutions
of Kolmogorov equation is usually considered the first step in proving
uniqueness for the corresponding stochastic equation (at least uniqueness in law, if not
pathwise uniqueness). However, even though the regularity of the solutions
constructed by~\cite{DaDeb2003JdMPeA} is high in terms of differentiability, the
regularity of their derivatives as functions on the infinite dimensional space
is not good enough, being defined only in subspaces where the weak solutions
of Navier-Stokes equations do not live continuously in time. Thus a careful
consideration of the assumptions does not allow to apply It\^{o} formula to
the composition of the solution of Kolmogorov equation and a weak solution of
the Navier-Stokes equations, a basic step in the usual proof of uniqueness. Said
differently, if the Kolmogorov equation is seen as the dual of the stochastic
equation (precisely the dual of the associated Fokker-Planck equation), the spaces
where the solutions of the two problems live, are not dual of one another, and thus any argument for uniqueness based on duality fails. In spite of this, the result of
\cite{DaDeb2003JdMPeA} can be considered the closest one to the solution of the
open problem. We know by personal communication that Giuseppe Da Prato always
kept in mind the open problem and continued to identify potential paths to its solution.

Below we describe some side results that may enrich the previous picture.
Inspired by the results in~\cite{DaDeb2003JdMPeA, DebOda2006Jee}, a theory of
Markov selections for stochastic 3D Navier-Stokes equations was developed in~\cite{FlaRom2008PTRF},
with a special property of continuous dependence on initial conditions that is
unique with respect to the deterministic case and thus worthy to be
mentioned; this is Section \ref{sect 3D NS} of this paper. Having
touched the difficulty to advance with additive noise in proving
regularization by noise, around 2010 there has been a shift to other kinds of
noise. Among them, multiplicative noise of transport type occupied a relevant
position (but it is not the only example; see for instance a noise multiplying
the Laplacian in Schr\"{o}dinger equation \cite{DebTsu2011JdMPeA} or multiplying the
nonlinear term in Hamilton Jacobi equations and conservation laws
\cite{GasGes2019PTRF, GesMau2018CPDE}). Heuristically, a multiplicative transport noise
is the Eulerian counterpart to an additive noise at the Lagrangian level, hence
could transfer the special well posedness properties of additive noise for
finite dimensional systems to the case of PDEs. The first results have been
for linear transport and advection equations \cite{FlaGubPri2010IM, FlaMauNek2014JMFM} but also
special solutions (point concentrated) of 2D Euler equations and 1D
Vlasov-Poisson system have been regularized by a similar noise \cite{FlaGubPri2011SPatA, DelFlaVin2014CPAM}. The same was proved for a Leray-$\alpha$ model
\cite{BarBesFer2014SPatA} and, finally, in~\cite{BarFlaMor2010PAMS, Bia2013ECP} for dyadic models of turbulence, a topic that we shall review in
Section \ref{sect dyadic}.

All the cited results of regularization by noise due to multiplicative noise
of transport type have been for inviscid problems and that seemed to be a
rule. However, it was recently understood that such noise may have a regularizing
effect also on viscous problems, in particular the 3D Navier-Stokes equations
since it increases dissipation \cite{FlaLuo2019}. This is described in Section
\ref{sect multiplicative 3D NS} and it is conceptually of interest also
because the regularization is not in the form of restored uniqueness by noise  - as all the previously
mentioned results - but in the form of suppression of blow-up, the second open question
mentioned above,~\cite{Fef2006TMPP}. We hope this picture may convince young
researchers that there is still space for improvements, although the research
on this topic is slow and rarely based on repeated schemes.

\section{The 3D Navier-Stokes equations}
\label{sect 3D NS}

\subsection{Deterministic case}
\label{sect deterministic}

For simplicity of exposition we assume that the fluid lives on the torus
$\mathbb{T}^{3}=\mathbb{R}^{3}/\mathbb{Z}^{3}$. We will denote by $H$ (resp.~$V$) the
Hilbert space of $L^{2}\left(  \mathbb{T}^{3},\mathbb{R}^{3}\right)  $ (resp. $W^{1,2}\left(  \mathbb{T}^{3},\mathbb{R}^{3}\right)  $) divergence free zero
average vector fields (see~\cite{Tem1995} for more precise details about the
boundary conditions). Let us recall, among others, the following basic results from~\cite{Tem1995}:
\begin{enumerate}
	\item Given $u_{0}\in H$, there exists a weak solution, namely a function of
	class%
	\[
	u\in L^{\infty}\left(  0,T;H\right)  \cap L^{2}\left(  0,T;V\right)
	\]
	weakly continuous in $H$, satisfying the identity%
	\begin{equation*}
	\left\langle u\left(  t\right)  ,\phi\right\rangle +\nu\int_{0}^{t}%
	\left\langle \nabla u\left(  s\right)  ,\nabla\phi\right\rangle
	\mathrm{d} s=\left\langle u_{0},\phi\right\rangle +\int_{0}^{t}\left\langle u\left(
	s\right)  ,u\left(  s\right)  \cdot\nabla\phi\right\rangle \mathrm{d} s 
	\end{equation*}
	for every $\phi\in V$;
	\item If $u_{0}\in V$, there exists a unique \textit{maximal} solution $u\in
	C\left(  [0,\tau);V\right)  $.
\end{enumerate}

Two questions (remember that we are in dimension $3$) remain open and represent fundamental problems in PDE
theory (see once more~\cite{Fef2006TMPP}):

\begin{enumerate}
\item Are weak solutions unique?

\item When $u_{0}\in V$, do we have $\tau=+\infty$ or%
\[
\tau<\infty\text{, }\lim_{t\uparrow\tau}\left\Vert u\left(  t\right)
\right\Vert _{V}=+\infty?
\]

\end{enumerate}

Here and in the following, we denote by $\left\Vert \cdot\right\Vert _{H}$ and $\left\Vert
\cdot\right\Vert _{V}$ the usual norms in $H$ and $V$ respectively and by
$\left\langle \cdot,\cdot\right\rangle $ either the scalar product in $H$ or
its extension to a dual pairing between spaces in duality with respect to $H$.

\subsection{Stochastic case, additive noise}

Generalization of the result of existence of weak solutions to the stochastic
case, with different types of noise, are now well-known, see for instance
\cite{FlaGat1995PTRF}, \cite{Fla2008CIME} and references therein. Let us mention some
elements in the case of additive noise. The formal notation is%
\begin{align*}
\mathrm{d} u+\left(  u\cdot\nabla u+\nabla p\right)  \mathrm{d} t &  =\Delta u\mathrm{d} t+\mathrm{d} W_{t}\\
\operatorname{div}u &  =0.
\end{align*}
Since space-time white noise is particularly attractive thanks to the outstanding
contributions of the theory of regularity structures and paracontrolled
distributions, let us first discuss this case, also because the general
results existing in the literature, for simpler nonlinearities, of
regularization by noise (like \cite{Gyo1998SPatA}, \cite{DaFla2010JoFA}) usually assume $W$ to
be a space-time white noise, namely a formal expression of the form%
\[
W_{t}\left(  x\right)  =\sum_{k\in\mathbb{Z}_{0}^{3},\alpha=1,2}\beta
_{t}^{k,\alpha}e_{k,\alpha}\left(  x\right)
\]
where the series converges in mean square in a distributional
space.
Here $\mathbb{Z}_{0}^{3}$ is $\mathbb{Z}^{3}\backslash\left\{  0\right\}  $ and
$\left(  e_{k,\alpha}\right)  _{k\in\mathbb{Z}_{0}^{3},\alpha=1,2}$ is a
complete orthonormal system of $H$ of the form%
\[
e_{k,\alpha}\left(  x\right)  =a_{k,\alpha}e^{2\pi ik\cdot x}\qquad
k\in\mathbb{Z}_{0}^{3},\alpha=1,2
\]
where $a_{k,1},a_{k,2}$ is an orthonormal basis of the plane perpendicular to
$k$ in $\mathbb{R}^{3}$. Finally, $\left(  \beta^{k,\alpha}\right)  _{k\in
\mathbb{Z}_{0}^{3},\alpha=1,2}$ is a family of complex Brownian motions
defined as follows: we take a family $\left(  W^{k,\alpha}\right)  _{k\in\mathbb{Z}%
_{0}^{3},\alpha=1,2}$ of real independent Brownian motions, we partition
$\mathbb{Z}_{0}^{3}$ in two sets $\mathbb{Z}_{+}^{3}$ and $\mathbb{Z}_{-}%
^{3}=-\mathbb{Z}_{+}^{3}$, and for all $k\in\mathbb{Z}_{+}^{3}$ we set
$\beta_{t}^{k,\alpha}=W^{k,\alpha}+iW^{-k,\alpha}$;\ for $k\in\mathbb{Z}%
_{-}^{3}$ we set $\beta_{t}^{k,\alpha}=W^{-k,\alpha}-iW^{k,\alpha}$. 

However, the solution of a parabolic equation
with space-time white noise is expected to be a function, and not just a
distribution, only when the spatial dimension is 1. In dimension 2 it is expected to be a distribution of Sobolev
class $H^{-\epsilon}$. This case was successfully investigated for
Navier-Stokes equations by Da Prato and Debussche in a seminal paper
\cite{DaDeb2003JdMPeA}; however, it is a 2D case, not competitive with the
deterministic theory (although striking from the stochastic viewpoint for
several reasons). In dimension $3$ the solutions are expected to be distributions of
cass $H^{\frac{1}{2}-\epsilon}$. A theorem of existence in such very singular
regime has been proven in~\cite{ZhuZhu2015JoDE}, but its relevance in view of a full
well posedness result is not clear. Thus we shall always consider more regular
noises, usually satisfying at least  the property that $W_{t}$ itself is a stochastic
process in $H$, an assumption achieved by requiring%
\begin{equation}
W_{t}\left(  x\right)  =\sum_{k\in\mathbb{Z}_{0}^{3},\alpha=1,2}%
\sigma_{k,\alpha}\beta_{t}^{k,\alpha}e_{k,\alpha}\left(  x\right)\;,
\label{noise}%
\end{equation}
where $\left(  \sigma_{k,\alpha}\right)  _{k\in\mathbb{Z}_{0}^{3},\alpha=1,2}$
are real numbers satisfying
\begin{equation}
\sum_{k\in\mathbb{Z}_{0}^{3},\alpha=1,2}\sigma_{k,\alpha}^{2}<\infty
.\label{trace class}%
\end{equation}

As said above, with this choice of noise, the equation was considered by several authors,
see for instance~\cite{FlaGat1995PTRF}. One can give a weak formulation as
\begin{equation*}
\left\langle u\left(  t\right)  ,\phi\right\rangle +\nu\int_{0}^{t}%
\left\langle \nabla u\left(  s\right)  ,\nabla\phi\right\rangle
\mathrm{d} s=\left\langle u_{0},\phi\right\rangle +\int_{0}^{t}\left\langle u\left(
s\right)  ,u\left(  s\right)  \cdot\nabla\phi\right\rangle \mathrm{d} s+\left\langle
W_{t},\phi\right\rangle\;, 
\end{equation*}
asking that
\begin{align*}
u  &  \in L_{\mathcal{F}}^{2}\left(  \Omega;\mathcal{H}\right) \\
\mathcal{H}  &  :=L^{2}\left(  0,T;W^{1,2}\right)  \cap C_{w}\left(
0,T;L^{2}\right) \;,
\end{align*}
namely that, on a probability space $\left(  \Omega,\mathcal{A},P\right)  $
with a filtration $\mathcal{F=}\left(  \mathcal{F}_{t}\right)  _{t\geq0}$ and
Brownian motions $\beta_{t}^{k}$ adapted to the filtration, $u$ is a weakly
continuous $\left(  \mathcal{F}_{t}\right)  $-adapted process in $H$, with
paths also of class $L^{2}\left(  0,T;V\right)  $, with suitable square
integrability properties (not needed here in detail), such that for all
$\phi\in V$ the previous identity holds true uniformly in time, with
probability one. When the tuple $\left(  \Omega,\mathcal{A},P,\left(
\mathcal{F}_{t}\right)  ,\left(  \beta_{t}^{k}\right)  \right)  $ is not
prescribed a priori, we say that a weak solution is \textit{a weak martingale
solution}. The existence of weak martingale solutions, as said above, is
known\footnote{When the tuple is arbitrarily given a priori, existence of
solutions is called strong existence; strong existence is open for the 3D
Navier-Stokes equations with additive noise.}.

Several extensions of more sophisticated deterministic results have been
proved in this stochastic setting. Among them, let us recall a generalization of the theory of Hausdorff
dimension of the set of singular points, the theory of \textit{
Caffarelli-Kohn-Nirenberg}. In the deterministic case, it claims that the set
$S$ of singular points in time-space may have at most Hausdorff dimension 1,
with 1-dimensional Hausdorff measure equal to zero. A full generalization to
the stochastic case has been obtained in~\cite{FlaRom2002TAMS}, with the following
probabilistic improvement.

\begin{theorem}
For stationary solutions (deterministic or stochastic case), if $S_{t}$ is the
random set of singularities at time $t$, then
\[
P\left(  S_{t}=\emptyset\right)  =1\;,
\]
for all $t\geq0$.
\end{theorem}

\subsubsection{Role of Kolmogorov equation for uniqueness in law}

For stochastic equations, uniqueness in law is the property stating that any two
solutions, possibly constructed on different probability spaces, have the same
law. This property is weaker than pathwise uniqueness, which is itself weaker than
path by path uniqueness At the same time, it is stronger than the uniqueness of the
associated Fokker-Planck equation. (We do not discuss these definitions here.)

How could one prove uniqueness in law by means of probabilistic arguments? Girsanov
theorem is the easiest method but it cannot work for Navier-Stokes equations,
as Ferrario has shown in~\cite{Fer2016EJP}. In general it seems that the
Girsanov approach has limitations that are too strong. The Kolmogorov approach, on the other hand, is more
flexible. The rough ``principle'' is that:

\begin{enumerate}[label=\roman*]
	\item If we control one derivative of Kolmogorov solution, we may try to prove
uniqueness in law;
	\item If we control two derivatives, and the first one is uniformly bounded, we
may try to prove pathwise uniqueness.
\end{enumerate}

Let us see some details on this topic. Consider an abstract stochastic
equation in Hilbert space:
\[
	\mathrm{d} u=\left(  Au+B\left(  u\right)  \right)  \mathrm{d} t+\mathrm{d} W_{t}%
\]
where $A$ is a negative selfadjoint operator, $B$ satisfies suitable assumptions, and $W$ is a
Brownian motion in $H$ with trace-class covariance $Q$ (the noise~\eqref{noise} is of this form when~\eqref{trace class} holds). Consider the
infinite dimensional backward Kolmogorov equation%
\begin{align*}
\partial_{t}U+\frac{1}{2}Tr\left(  QD^{2}U\right)  +\left\langle Au+B\left(
u\right)  ,DU\right\rangle  &  =0\\
U|_{t=T} &  =\phi\;,
\end{align*}
where for the time being we do not give precise definitions of the
single objects\footnote{Just notice, with a certain degree of formality, that
$U=U\left(  t,u\right)  $ is a real function defined on $\left[  0,T\right]
\times H$, with the notation $u\in H$; $DU\left(  t,u\right)  $ is its
differential in the $H$-variable, element of $H$, $\left\langle Au+B\left(
u\right)  ,DU\left(  t,u\right)  \right\rangle $ is its scalar product in $H$
with the vector $Au+B\left(  u\right)  $, $D^{2}U\left(  t,u\right)  $ is the
second differential, an operator on $H$, and $Tr\left(  QD^{2}U\left(
t,u\right)  \right)  $ is the trace of the operator $QD^{2}U\left(
t,u\right)  $; finally $\phi$ is a real function on $H$}. Heuristically,
assume that the Kolmogorov equation has a sufficiently smooth solution and assume that
$u\left(  t\right)  $ is a solution of the stochastic equation. By It\^{o}
formula, for $0\leq r\leq t\leq T$,%
\begin{multline*}
U\left(  t,u\left(  t\right)  \right)  -U\left(  r,u\left(  r\right)  \right)
  =\int_{0}^{T}\left\langle DU\left(  s,u_{s}\right)  ,\mathrm{d} W_{s}\right\rangle \\
  +\int_{r}^{t}\left(  \partial_{s}U+\frac{1}{2}Tr\left(  QD^{2}U\right)
+\left\langle A\cdot+B\left(  \cdot\right)  ,DU\right\rangle \right)  \left(
s,u\left(  s\right)  \right)  \mathrm{d} s\;,
\end{multline*}
and thus, by the Kolmogorov equation,%
\[
U\left(  t,u\left(  t\right)  \right)  -U\left(  r,u\left(  r\right)  \right)
=\int_{0}^{T}\left\langle DU\left(  s,u_{s}\right)  ,dW_{s}\right\rangle\;.
\]
If $DU$ is good enough to have
\begin{equation}
\mathbb{E}\int_{0}^{T}\left\Vert DU\left(  s,u_{s}\right)  \right\Vert
_{H}^{2}ds<\infty\;, \label{integrability}%
\end{equation}
then $\mathbb{E}\int_{0}^{T}\left\langle DU\left(  s,u_{s}\right)
,dW_{s}\right\rangle =0$ and we deduce%
\[
\mathbb{E}\phi\left(  u\left(  T\right)  \right)  =\mathbb{E}U\left(
0,u_{0}\right)  .
\]
This, by the arbitrariness of $\phi$, identifies the law of $u\left(
T\right)  $ (and $T$ is arbitrary). With more work, as explained for instance
in \cite{StrVar2007}, we identify the law of the process. Let us remark that
Giuseppe Da Prato was the main investigator of Kolmogorov equations in
infinite dimensional spaces, see for instance his two books \cite{Da2012, DaZab2002}.

The classical idea to investigate Kolmogorov equations in infinite dimensions
is by perturbation. In order to describe it, let us reverse time by
setting $V\left(  t\right)  =U\left(  T-t\right)  $; now we have to study
the forward equation%
\begin{align*}
\partial_{t}V &  =\frac{1}{2}Tr\left(  QD^{2}V\right)  +\left\langle
Au+B\left(  u\right)  ,DV\right\rangle \\
V|_{t=0} &  =\phi.
\end{align*}
Introducing the Gaussian semigroup solving
\[
\partial_{t}\mathcal{S}_{t}\phi=\frac{1}{2}Tr\left(  QD^{2}\mathcal{S}_{t}%
\phi\right)  +\left\langle Au,D\mathcal{S}_{t}\phi\right\rangle\;, 
\]
with $\mathcal{S}_{0}\phi=\phi$,
one rewrites the equation in perturbative form%
\[
V\left(  t\right)  =\mathcal{S}_{t}\phi+\int_{0}^{t}\mathcal{S}_{t-s}%
\left\langle B\left(  u\right)  ,DV\left(  s\right)  \right\rangle \mathrm{d} s.
\]
In order to apply a fixed point argument to this equation in suitable spaces, 
it is necessary to have good gradient bounds on the Gaussian semigroup. Those
usually proved, under suitable assumptions on the pair $\left(  A,Q\right)  $,
have the form%
\begin{equation}
\left\Vert D\mathcal{S}_{t}\phi\right\Vert _{0}\leq\frac{C}{t^{\gamma}%
}\left\Vert \phi\right\Vert _{0}\;, \label{eq:grad_bound}
\end{equation}
with $\gamma\in\left(  0,1\right)  $. Here $\left\Vert \phi\right\Vert _{0}$
is the uniform norm of a function or a vector defined on $H$. Unfortunately, a
great limitation of this perturbative approach is that $B$ has to be bounded,
see for instance \cite{DaFla2010JoFA, DaFlaPriRoc2013AP, DaFlaRocVer2016AP}. Moreover, the assumptions on
$\left(  A,Q\right)  $ to have the gradient bound~\eqref{eq:grad_bound} are far from those satisfied
by the linear part of 3D Navier-Stokes equations.

Da Prato and Debussche in~\cite{DaDeb2003JdMPeA} made a breakthrough on this topic in
the direction of 3D Navier-Stokes equations: under suitable assumptions on the
coefficients $\sigma_{k}$ (the idea behind the assumptions is that the coefficients cannot go
to zero too fast), they discovered a way to construct smooth solutions of the
associated infinite dimensional Kolmogorov equation. Without pretending to
explain in a sentence the very elaborate procedure developed in
\cite{DaDeb2003JdMPeA}, let us only mention that it starts with the very innovative
idea of introducing a penalized evolution operator $\mathcal{R}\left(
s,t\right)  $ in place of the Gaussian semigroup:
\[
U\left(  t\right)  =\mathcal{R}\left(  0,t\right)  \phi+\int_{0}%
^{t}\mathcal{R}\left(  s,t\right)  \left(  \left\langle B\left(  u\right)
,DU\left(  s\right)  \right\rangle -V\left(  s\right)  \right)  \mathrm{d} s\;.
\]
Using this method it is possible to prove the existence of a smooth solution
$U\left(  t,u\right) $.

The solution $U$ is differentiable (in fact twice differentiable), but with bounds
on derivatives of the form%
\[
\left\langle h,DU\left(  t,u\right)  \right\rangle \leq C\left(  t\right)
\left\Vert h\right\Vert _{W^{2,2}}\left(  1+\left\Vert u\right\Vert _{W^{2,2}%
}\right)\;,
\]
namely depending on a Sobolev norm in the infinite dimensional variable $u$,
which is quite demanding from the viewpoint of the regularity of solutions of
3D Navier-Stokes equations. If we go back to the sufficient condition~\eqref{integrability}, we see that weak solutions do not have sufficient
regularity. In principle there could be several weaker ways to proceed, which
do not require directly~\eqref{integrability}, but no way has been found yet.

Technically, Da Prato-Debussche \cite{DaDeb2003JdMPeA} is one of the most advanced
works on stochastic 3D Navier-Stokes equations. Not only it constructs
solutions to the Kolmogorov equation, but it also identifies two new
properties: Markov selections and strong Feller property, discussed below.

\subsection{Small times versus large times}

Let $u_{0}\in L^{2}$ be an initial condition and $u\in L_{\mathcal{F}}%
^{2}\left(  \Omega;\mathcal{H}\right)  $ be a (possibly non-unique) weak
solution. Using the properties of conditional expectation, let us decompose%
\begin{align*}
\mathbb{E}\left[  \phi\left(  u\left(  T\right)  \right)  \right]   &
=\mathbb{E}\left[  \mathbb{E}\left[  \phi\left(  u\left(  T\right)  \right)
|u\left(  t_{0}\right)  \right]  \right]  \\
&  =\int_{H}\mathbb{E}\left[  \phi\left(  u\left(  T\right)  \right)
|u\left(  t_{0}\right)  =v\right]  \mu_{t_{0}}\left(  \mathrm{d} v\right)\;,
\end{align*}
where $\mu_{t_{0}}$ is the law of $u\left(  t_{0}\right)  $, and $\phi$ is a
smooth test functional on $H$. The hope is to propagate good properties, which
hold for small times, to large times.

Assume for every initial condition $u_{0}$ we \textit{select} a weak solution
$u\left(  \cdot;u_{0}\right)  \in L_{\mathcal{F}}^{2}\left(  \Omega
;\mathcal{H}\right)  $. Uniqueness is not known, but we may make selections,
following different criteria; the simplest one is measurable-in-$u_{0}$
selection, but a more refined one, following \cite{StrVar2007}, is a Markov
selection (see below). For each one of the selected solutions we have the
decomposition above%
\begin{equation}
\mathbb{E}\left[  \phi\left(  u\left(  T;u_{0}\right)  \right)  \right]
=\int_{L^{2}}\mathbb{E}\left[  \phi\left(  u\left(  T;u_{0}\right)  \right)
|u\left(  t_{0};u_{0}\right)  =v\right]  \mu_{t_{0},u_{0}}\left(  \mathrm{d} v\right)\;,
\label{decomposition}%
\end{equation}
where $\mu_{t_{0},u_{0}}$ is the law of $u\left(  t_{0};u_{0}\right)  $.
One can already notice the germ of a special property: if $u_{0}\in V$ and $t_{0}$ is
small enough, the law $\mu_{t_{0},u_{0}}\left(  \mathrm{d} v\right)  $ is
`` almost'' independent of the selection,
since for $u_{0}\in V$ the solution is locally unique. The limitation
`` almost'' refers to the fact that
`` locally'', in the stochastic case, means
randomly local, hence we know uniqueness up to time $t_{0}$ only with large probability.

Assume $u_{0}^{n},u_{0}\in V$ are such that
\[
u_{0}^{n}\overset{V}{\rightarrow}u_{0}.
\]
In the deterministic case, one can find $t_{0}$ small enough that unique
solutions $u^{n},u$ exist on $\left[  0,t_{0}\right]  $ with initial
conditions $u_{0}^{n},u_{0}$ and $u^{n}\rightarrow u$ in $C\left(  \left[
0,t_{0}\right]  ;V\right)  $. In the stochastic case, a similar result holds
with large probability~\cite{FlaRom2008PTRF}: for every $\epsilon>0$ there
exists $t_{0}>0$ such that solutions exist, and are pathwise unique, in $C\left(
\left[  0,t_{0}\right]  ;V\right)  $ with probability greater than
$1-\epsilon$; at the same time, $u^{n}\rightarrow u$ in $C\left(  \left[  0,t_{0}\right]
;V\right)  $ with probability greater than $1-\epsilon$. Forgetting about this
$\epsilon$ for the sake of simplicity of the heuristic explanation (the details
are in~\cite{FlaRom2008PTRF}), we have
\[
u\left(  t_{0};u_{0}^{n}\right)  \underset{\text{a.s.}%
}{\overset{V}{\rightarrow}}u\left(  t_{0};u_{0}\right)\;, 
\]%
and 
\begin{equation}
\lim_{u_{0}^{n}\rightarrow u_{0}}\int\psi\left(  v\right)  \mu_{t_{0}%
,u_{0}^{n}}\left(  \mathrm{d} v\right)  =\int\psi\left(  v\right)  \mu_{t_{0},u_{0}%
}\left(  \mathrm{d} v\right)  \label{weak convergence}%
\end{equation}
for a large class of continuous functions $\psi$. 

The previous result is only the stochastic analog of a deterministic property
of local well posedness. But in the stochastic case it is here that we have
more. Under strong assumptions on the noise (the same ones that allowed to solve 
the Kolmogorov equation in~\cite{DaDeb2003JdMPeA}), \textit{strong Feller
property} holds at time $t_{0}$ (again we simplify the exposition forgetting
about a small probability $\epsilon$ of having a different property)%
\[
\lim_{u_{0}^{n}\rightarrow u_{0}}\mu_{t_{0},u_{0}^{n}}=\mu_{t_{0},u_{0}}\text{
in total variation.}%
\]
Convergence in total variation essentially means that~\eqref{weak convergence}
is extended to a large class of measurable functions, something impossible in
the deterministic case, where $\mu_{t_{0},u_{0}^{n}}$ and $\mu_{t_{0},u_{0}}$
are delta Dirac masses! Using the decomposition property~\eqref{decomposition}, 
one can prove:

\begin{theorem}
Assume $\mathbb{E}\left[  \phi\left(  u\left(  T;u_{0}\right)  \right)
|u\left(  t_{0};u_{0}\right)  =v\right]  $ is independent of $u_{0}$ and there
exists a function $g_{\phi}\left(  T,t_{0},v\right)  $, measurable in $v$,
such that
\[
g_{\phi}\left(  T,t_{0},v\right)  =\mathbb{E}\left[  \phi\left(  u\left(
T;u_{0}\right)  \right)  |u\left(  t_{0};u_{0}\right)  =v\right]  .
\]
Then $\mathbb{E}\left[  \phi\left(  u\left(  T;u_{0}^{n}\right)  \right)
\right]  \rightarrow\mathbb{E}\left[  \phi\left(  u\left(  T;u_{0}\right)
\right)  \right]  $, namely continuous dependence propagates to large times.
\end{theorem}

The assumption of the theorem, existence of $g_{\phi}\left(  T,t_{0},v\right)
$, is essentially the Markov property. The question is: can we make a
selection which satisfies the Markov property?

Yes, following~\cite{DaDeb2003JdMPeA, DebOda2006Jee, FlaRom2008PTRF} we know:

\begin{theorem}
For 3D Navier-Stokes, Markov selections exist. If the noise is strong enough,
they are strong Feller, hence solutions depend continuously on the initial
conditions, also for large times, in the topology of $V$.
\end{theorem}

The previous theorem can be considered the most advanced innovative result of
the stochastic theory with respect to the deterministic one. Nothing like this
theorem is known in the deterministic case. 

Can we do more? The following trick in semigroup theory is well known: if
$A:D\left(  A\right)  \subset H\rightarrow H$ generates a strongly continuous
semigroup $S_{t}$, $t\geq0$, and $u\left(  t\right)  $ solves $u^{\prime
}\left(  t\right)  =Au\left(  t\right)  $, then
\[
u\left(  t\right)  =S_{t}u\left(  0\right)  .
\]
Indeed,%
\[
\frac{\mathrm{d}}{\mathrm{d} s}S_{t-s}u\left(  s\right)  =-AS_{t-s}u\left(  s\right)
+S_{t-s}Au\left(  s\right)  =0.
\]
In other words: when we have a strongly continuous flow, all solutions
coincide with those of the flow. Such uniqueness result, however, holds in the
framework of semigroup theory; it is only heuristically a general principle.
In the case described above, we have something similar concerning the
assumptions: we have a Markov, strong Feller, selection. But, in spite of
many attempts, we have not found a rigorous way to deduce that it
``incorporates'' every weak solution. 

The Markov strong Feller selection is a priori not unique and, based on
results proved in~\cite{StrVar2007} in an easier context than the Navier-Stokes
equations, we should expect uniqueness of Markov selections if and only if
there is uniqueness of individual solutions. It is however possible that some
Markov selection carries more specific information and may be elevated to a
special role. Sufficient conditions for uniqueness of Markov selections are
given in~\cite{FlaRom2008PTRF, Rom2008JSP}.

\subsection{Multiplicative transport noise}
\label{sect multiplicative 3D NS}

Another noise received increasing attention in fluid mechanics problems. It is
inspired by the transport term $u\cdot\nabla u$ and has the form (compare with~\eqref{noise})%
\[
\nabla u\circ\mathrm{d} W=\sum_{k\in\mathbb{Z}_{0}^{3},\alpha=1,2}\sigma_{k,\alpha
}\left(  e_{k,\alpha}\cdot\nabla u\right)  \circ \mathrm{d}\beta_{t}^{k,\alpha}.
\]
The multiplication is understood in the Stratonovich sense, recognized to be
the right one throughout the literature on this subject (e.g.~\cite{Hol2015PotRSAMPaES, MolRuzSok1985SPU, MajTimEij2001CPAM}). A short introduction to this detail can be
found in \cite{FlaMauNek2014JMFM}.

In a sense, the velocity field $u$ which transports other quantities (like $u$
itself in $u\cdot\nabla u$, or terms like $u\cdot\nabla T$ in heat transport)
is replaced by $u+W$. The resulting stochastic Navier-Stokes equations are
\begin{align*}
\mathrm{d} u+\left(  u\cdot\nabla u+\nabla p\right)  \mathrm{d} t &  =\Delta u \mathrm{d} t+\sum
_{k\in\mathbb{Z}_{0}^{3},\alpha=1,2}\sigma_{k,\alpha}\left(  e_{k,\alpha}%
\cdot\nabla u\right)  \circ \mathrm{d}\beta_{t}^{k,\alpha}\\
\operatorname{div}u &  =0.
\end{align*}
There is another, non-equivalent, way to introduce transport noise; it is at
the level of the equation for the vorticity $\xi=\operatorname{curl}u$, which
in the case of 3D deterministic Navier-Stokes equations is%
\[
\partial_{t}\xi+u\cdot\nabla\xi=\Delta\xi+\xi\cdot\nabla u\;,
\]
also written, using the Lie derivative $\mathcal{L}_{u}\xi=u\cdot\nabla\xi
-\xi\cdot\nabla u$, as%
\[
\partial_{t}\xi+\mathcal{L}_{u}\xi=\Delta\xi.
\]
The natural perturbation of this equation is%
\begin{align*}
\mathrm{d}\xi+u\cdot\nabla\xi \mathrm{d} t &  =\Delta\xi \mathrm{d} t+\xi\cdot\nabla u\mathrm{d} t+\sum
_{k\in\mathbb{Z}_{0}^{3},\alpha=1,2}\sigma_{k,\alpha}\left(  e_{k,\alpha}%
\cdot\nabla\xi\right)  \circ \mathrm{d}\beta_{t}^{k,\alpha}\\
&  -\sum_{k\in\mathbb{Z}_{0}^{3},\alpha=1,2}\sigma_{k,\alpha}\left(  \xi
\cdot\nabla e_{k,\alpha}\right)  \circ \mathrm{d}\beta_{t}^{k,\alpha}%
\end{align*}
considered in \cite{Hol2015PotRSAMPaES, CriFlaHol2019JNS}: it corresponds to the replacement of $u$
with $u+W$ in the Lie derivative (which corresponds to the same replacement at
the Lagrangian level)
\begin{align*}
\mathcal{L}_{u}\xi \mathrm{d} t  & \rightarrow\mathcal{L}_{u}\xi \mathrm{d} t+\mathcal{L}_{\circ
\mathrm{d} W}\xi\\
& :=\mathcal{L}_{u}\xi \mathrm{d} t+\sum_{k\in\mathbb{Z}_{0}^{3},\alpha=1,2}%
\mathcal{L}_{\sigma_{k,\alpha}e_{k,\alpha}\circ \mathrm{d}\beta_{tt}^{k,\alpha k}}\xi.
\end{align*}
The nonlinearity is composed, at the vorticity level, of two terms: the
transport of vorticity $u\cdot\nabla\xi$ and the vortex stretching $\xi
\cdot\nabla u$. Accordingly, in the previous equation there is an additional
stochastic transport and stochastic stretching. When vorticity is replaced by
magnetic moment, this stochastic perturbation was considered in the framework
of the dynamo theory in~\cite{MolRuzSok1985SPU}. Notice that in the 2D case stretching
cannot occur, since the vorticity is orthogonal to the plane of fluid motion,
hence the equation reduces to (see for instance \cite{BrzFlaMau2016ARMA})%
\[
\mathrm{d}\xi+u\cdot\nabla\xi \mathrm{d} t=\Delta\xi \mathrm{d} t+\sum_{k\in\mathbb{Z}_{0}^{3}}\sigma
_{k}\left(  e_{k}\cdot\nabla\xi\right)  \circ \mathrm{d}\beta_{t}^{k}
\]
It is worth noticing that it is not necessary anymore to sum over the index $\alpha=1,2$ because the
linear space orthogonal to $k$ is now a line). Concerning motivations for the
model with transport noise, let us mention model reduction, see \cite{MajTimEij2001CPAM}, in addition
to other motivations like \cite{BrzCapFla1992SAA}, \cite{MikRoz2005AP}, and
more recently \cite{Hol2015PotRSAMPaES}.

Starting from 2010, several simpler models proved to be regularized by
transport noise, as already remarked in the Introduction: linear transport and
advection equations \cite{FlaGubPri2010IM, FlaMauNek2014JMFM, FlaOli2018JEE}, special
solutions of 2D Euler equations and 1D Vlasov-Poisson system \cite{FlaGubPri2011SPatA, DelFlaVin2014CPAM}, Leray-$\alpha$ model \cite{BarBesFer2014SPatA} and, as more extensively
discussed below in Section \ref{sect dyadic}, dyadic models of turbulence
\cite{BarFlaMor2010PAMS}, \cite{Bia2013ECP}. In all these cases the PDE is
inviscid. But it was recently understood that such a noise may have a
regularizing effect also on viscous problems, in particular the 3D
Navier-Stokes equations because it increases dissipation. Let us
briefly summarize this result, from~\cite{FlaLuo2019}.

The first important remark is that it holds for a sort of artificial
modification of the noise above: we consider only the stochastic transport
term - as in the 2D case -, neglecting the stochastic stretching term, but
maintaining the 3-dimensionality of the equation. The precise model is%
\begin{equation}
\mathrm{d}\xi+u\cdot\nabla\xi \mathrm{d} t=\Delta\xi \mathrm{d} t+\xi\cdot\nabla u \mathrm{d} t+\sum_{\substack{k\in
\mathbb{Z}_{0}^{3}\\ \alpha=1,2}}\sigma_{k,\alpha}\Pi\left(  e_{k,\alpha}%
\cdot\nabla\xi\right)  \circ \mathrm{d}\beta_{t}^{k,\alpha}
\label{3D NS transport noise Strat}%
\end{equation}
where $\Pi$ is the projection on divergence free fields, necessary since the
sum of all other terms is divergence free (notice that, on the contrary, the full noise
$\mathcal{L}_{\circ dW}\xi$ does not require projection since it is already
divergence free). In~\cite{FlaLuo2019} there is an attempt to motivate this choice of noise, 
but it remains true that the full noise $\mathcal{L}_{\circ dW}\xi$ is much more
natural, while at the same time the latter spoils the result of regularization by noise,
as shown in~\cite{FlaLuo2019}. This discrepancy will be the object of future investigation.

In order to understand the result in~\cite{FlaLuo2019}, let us recall the second open
problem presented in Section~\ref{sect deterministic}, restated here as follows:
when $\xi_{0}\in H$, do we have%
\[
\tau<\infty\text{, }\lim_{t\uparrow\tau}\left\Vert \xi\left(  t\right)
\right\Vert _{H}=+\infty?
\]
We have discovered that transport noise may improve the control of $\left\Vert
\xi\left(  t\right)  \right\Vert _{H}$. In the deterministic case, the norm
$\left\Vert \xi\left(  t\right)  \right\Vert _{H}^{2}$ can be controlled
\textit{locally} from
\[
\partial_{t}\xi+u\cdot\nabla\xi-\xi\cdot\nabla u=\Delta\xi\;,
\]
by energy type estimates:%
\[
\frac{1}{2}\frac{\mathrm{d}}{\mathrm{d} t}\left\Vert \xi\left(  t\right)  \right\Vert _{H}%
^{2}+\left\Vert \nabla\xi\left(  t\right)  \right\Vert _{H}^{2}=\left\langle
\xi\cdot\nabla u,\xi\right\rangle .
\]
The term $\left\langle \xi\cdot\nabla u,\xi\right\rangle $ describes the
\textit{stretching} of vorticity $\xi$ produced by the deformation tensor
$\nabla u$. This is the potential source of unboundedness of $\left\Vert
\xi\left(  t\right)  \right\Vert _{H}^{2}$. Sobolev and interpolation
inequalities give us (up to constants):%
\[
\left\langle \xi\cdot\nabla u,\xi\right\rangle \leq\left\Vert \xi\right\Vert
_{L^{3}}^{3}\leq\left\Vert \xi\right\Vert _{W^{1/2,2}}^{3}\leq\left\Vert
\xi\right\Vert _{L^{2}}^{3/2}\left\Vert \xi\right\Vert _{W^{1,2}}^{3/2}%
\leq\left\Vert \xi\right\Vert _{W^{1,2}}^{2}+\left\Vert \xi\right\Vert
_{L^{2}}^{6}%
\]
and this leads to%
\[
\frac{\mathrm{d}}{\mathrm{d} t}\left\Vert \xi\left(  t\right)  \right\Vert _{H}^{2}\leq
C\left\Vert \xi\right\Vert _{H}^{6}%
\]
which provides only a local control. However the interval of existence depends
on the \textit{viscosity coefficient} $\nu$: if we consider
\[
\partial_{t}\xi+u\cdot\nabla\xi-\xi\cdot\nabla u=\nu\Delta\xi\;,
\]%
the energy estimate become
\begin{align*}
\frac{1}{2}\frac{\mathrm{d}}{\mathrm{d} t}\left\Vert \xi\left(  t\right)  \right\Vert _{H}%
^{2}+\nu\left\Vert \nabla\xi\left(  t\right)  \right\Vert _{H}^{2} &
=\left\langle \xi\cdot\nabla u,\xi\right\rangle \\
&  \leq\left\Vert \xi\right\Vert _{L^{2}}^{3/2}\left\Vert \xi\right\Vert
_{W^{1,2}}^{3/2}\\
&  \leq\nu\left\Vert \nabla\xi\left(  t\right)  \right\Vert _{H}^{2}+\frac
{C}{\nu^{3}}\left\Vert \xi\right\Vert _{H}^{6}
\end{align*}%
leading in this case to
\[
\frac{\mathrm{d}}{\mathrm{d} t}\left\Vert \xi\left(  t\right)  \right\Vert _{H}^{2}\leq\frac
{C}{\nu^{3}}\left\Vert \xi\right\Vert _{H}^{6}.
\]
The explosion is delayed for large $\nu$. Not only that: beyond a threshold the
solution is global. This is the key for a regularization by noise: transport
noise improves dissipation, hence it delays blow-up.

Let us rewrite equation~\eqref{3D NS transport noise Strat} In It\^{o} form
(see~\cite{FlaMauNek2014JMFM} for an easy introduction to this operation):
\begin{align*}
\mathrm{d}\xi+u\cdot\nabla\xi \mathrm{d} t  & =\Delta\xi \mathrm{d} t+\xi\cdot\nabla u \mathrm{d} t+\sum
_{k\in\mathbb{Z}_{0}^{3},\alpha=1,2}\sigma_{k,\alpha}\Pi\left(  e_{k,\alpha
}\cdot\nabla\xi\right)  \mathrm{d}\beta_{t}^{k,\alpha}\\
& +\frac{1}{2}\sum_{k\in\mathbb{Z}_{0}^{3},\alpha=1,2}\sigma_{k,\alpha}^{2}%
\Pi\left(  e_{k,\alpha}\cdot\nabla\Pi\left(  e_{k,\alpha}\cdot\nabla
\xi\right)  \right)  \mathrm{d} t\;,
\end{align*}
where the stochastic term is now understood in It\^{o} sense. The corrector is
a pseudo-differential operator of second order, quite complicated
algebraically by the presence of the projector $\Pi$. Under suitable technical
conditions on the family of coefficients $\sigma=\left(  \sigma_{k,\alpha
}\right)  _{k\in\mathbb{Z}_{0}^{3},\alpha=1,2}$ (still quite general), the
corrector turns out to be of the form%
\[
\frac{1}{2}\sum_{k\in\mathbb{Z}_{0}^{3},\alpha=1,2}\sigma_{k,\alpha}^{2}%
\Pi\left(  e_{k,\alpha}\cdot\nabla\Pi\left(  e_{k,\alpha}\cdot\nabla
\xi\right)  \right)  =\nu_{\sigma}\Delta\xi+R_{\sigma}\left(  \xi\right)\;,
\]
where $\nu_{\sigma}>0$ is a coefficient depending on $\sigma$ and $R_{\sigma
}\left(  \xi\right)  $ is a quite complicated non-local second order
differential operator. The decomposition of the RHS as $\nu_{\sigma}\Delta\xi+R_{\sigma
}\left(  \xi\right)  $ is not purely artificial: the same corrector without
the two projections $\Pi$ would be simply equal to $\nu_{\sigma}\Delta\xi$;
 the remainder $R_{\sigma}\left(  \xi\right)  $ is what is left due to
the presence of the projections.

Now the key point is to parametrize $\sigma$ by a scaling parameter $N$:%
\[
\sigma^{N}=\left(  \sigma_{k,\alpha}^{N}\right)  _{k\in\mathbb{Z}_{0}%
^{3},\alpha=1,2}%
\]
in such a way that the corresponding coefficient $\nu_{\sigma^{N}}$ is independent
of $N$%
\[
\nu_{\sigma^{N}}=\nu
\]
and (this is the most difficult technical part of the work~\cite{FlaLuo2019})%
\[
\lim_{N\rightarrow\infty}R_{\sigma^{N}}\left(  \xi\right)  =-\frac{2}{5}%
\nu\Delta\xi.
\]
The solutions $\xi^{N}$ of the corresponding equation%

\begin{equation}
\mathrm{d}\xi^{N}+u^{N}\cdot\nabla\xi^{N}\mathrm{d} t=\Delta\xi^{N}\mathrm{d} t+\xi^{N}\cdot\nabla
u^{N} \mathrm{d} t+\sum_{k\in\mathbb{Z}_{0}^{3},\alpha=1,2}\sigma_{k,\alpha}^{N}%
\Pi\left(  e_{k,\alpha}\cdot\nabla\xi^{N}\right)  \circ \mathrm{d}\beta_{t}^{k,\alpha
}\label{3D NS transport noise Strat N}%
\end{equation}
will have the following properties, which are the main results of \cite{FlaLuo2019}.

\begin{theorem}
Let $\xi_{0}\in H$ and $\left[  0,T\right]  $ be given. In a suitable scaling
limit $N\rightarrow\infty$ corresponding to a sequence $\sigma^{N}$, $\xi^{N}$
converges in probability to the solution of%
\[
\partial_{t}\xi+\mathcal{L}_{u}\xi=\left(  1+\frac{5}{3}\nu\right)  \Delta\xi.
\]
It follows that for large $N$ the norm $\left\Vert \xi^{N}\left(  t\right)
\right\Vert _{H}^{2}$ is bounded on $\left[  0,T\right]  $, with high
probability (implying well posedness of $\xi^{N}$).
\end{theorem}

\begin{theorem}
Given $R_{0},\epsilon>0$, there exists $N$ with the following property: for
every initial condition $\xi_{0}\in H$ with $\left\Vert \xi_{0}\right\Vert
_{H}\leq R_{0}$, the stochastic 3D Navier-Stokes equations
(\ref{3D NS transport noise Strat N}) have a global unique solution, up to
probability $\epsilon$.
\end{theorem}

This result is a regularization by noise result because the viscosity in
equation~\eqref{3D NS transport noise Strat N} is 1 and, as discussed above
for the deterministic equations, with such viscosity only very small initial
conditions lead to global existence.

The previous results are inspired by several sources, among which we quote
\cite{ArnCraWih1983SJCO, BabMahNic1996EJMB, BarBesFer2014SPatA, ConKisRyzZla2008AM, FlaLuo2018, FlaGalLuo2019, Gal2019}.

\section{Regularization by noise in dyadic models}
\label{sect dyadic}

Even though the regularization by noise techniques did not work for 3D
Navier-Stokes equation, there are other equations that proved to be more
accessible with this tool. One special case, still in the area of
fluid-dynamics, is that of the dyadic models of turbulence.

\subsection{Dyadic models}


Shell models were introduced by the Russian school in the 1970s, as a
theoretical and computational tool to study the cascade phenomenon in
turbulent fluid dynamics. \ This is a mechanism (not yet completely
understood) that moves energy from one lengthscale to another, thus sustaining
turbulence. Richardson's cascade, also called direct energy cascade, moves the
energy from larger scales to smaller ones, whereas the inverse cascade moves
energy from smaller to larger scales, and seems to appear only in 2D turbulence.

The phenomenological idea behind the tree model proposed by Katz and
Pavlovi{\'c}~{\cite{KatPav2005TAMS}}, and called KP model in the next pages,
is the following: larger eddies in the turbulent fluid split into smaller ones
because of dynamical instabilities, and the kinetic energy moves from the
larger scales to the smaller ones. We simplify the picture by assuming that
eddies appear only at certain discrete scales, each the half of the previous
one. We also assume that the eddies fill the space, so that each eddy contains
$2^{d}$ eddies of the next scale.

In this way we have a tree structure, where each node is an eddy. Following
the notation introduced in~{\cite{BarBiaFlaMor2013JMP}}, if we denote by $J$
the set of nodes, each node $j$ has a set of children $\mathcal{O}_{j}$,
representing the smaller eddies generated by instability from~$j$. We call
generations the discrete scales where the eddies are, and denote the
generation of an eddy $j$ by $| j |$. At level (or generation) 0 we have the
single largest eddy, denoted with $\emptyset$, at generation 1 the $2^{d}$
eddies generated by the eddy at level 0 and so on. Also, we denote the parent
of a node $j$ by $\overline{\jmath}$.

Every node has a scalar quantity $X_{j}$ attached to it, the intensity of the
velocity field, with the square of this intensity being the kinetic energy. In
other words, the energy is the square of the $l^{2}$-norm: $\mathcal{E} (t) :=
\sum_{j} X_{j}^{2} (t)$. The intensities are coupled by the following
differential rules:
\begin{equation}
\label{eq:tree_dyadic}\dot{X}_{j} = - \nu\tilde{c}_{j} X_{j} + c_{j}
X_{\bar{\jmath}}^{2} - X_{j} \sum_{k \in\mathcal{O}_{j}} c_{k} X_{k},
\end{equation}
where we consider the coefficients $c_{j} = d_{j} 2^{\alpha| j |}$, with
$\alpha> 0$ and $d_{j} > 0$ for all $j \in J$ (and similarly for the
$\tilde{c}_{j} = \tilde{d}_{j} 2^{\gamma| j |}$), $d_{\emptyset} = 1$, and
$X_{\overline{\emptyset}} (t) \equiv f$, that is the forcing acts only on the
largest eddy. Most results are independent of the choice of $\alpha$, however
there are heuristic arguments that suggest $\alpha= \frac{d}{2} + 1$, which is
the value usually considered in the literature (see for
example~{\cite{KatPav2005TAMS,BarBiaFlaMor2013JMP}}).
In~{\cite{BarMorRom2014AP}} it was proven that $\alpha\leqslant\frac{5}{2}$
for a Littlewood-Paley decomposition of 3D Euler dynamics.

In~{\cite{KatPav2005TAMS}} and in~{\cite{BarBiaFlaMor2013JMP}}, $d_{j} = 1$
for all $j \in J$, but in~{\cite{BiaMor2017CMP}}, restricted to the inviscid
case (i.e.~$\nu= 0$) the coefficients $d_{j}$ are allowed to vary for
different nodes, with the assumption that $| \log d_{j} |$ is bounded. Moreover
a particular choice is introduced, the repeated coefficients models (or RCM),
in which the same fixed $2^{d}$ coefficients $\delta_{\omega}$ appear in every
set of siblings $\{ d_{k} : k \in\mathcal{O}_{j} \}$. For the RCM it is
possible to state and prove more interesting and deep results, due to its
simpler form.

Heuristically, we can think of this model as a (simplified) wavelet
decomposition of Navier-Stokes equations, see for
example~{\cite{KatPav2005TAMS,Che2008TAMS,BiaMor2017CMP}}. However, this is
not a rigorous derivation, as pointed out in~{\cite{Wal2006PAMS}}: the KP
model is constructed in such a way that it mimics Navier-Stokes (in particular
with respect to the energy cascade phenomenon).

If we choose to have only one intensity per shell, that is we consider all
nodes in a generation as collapsed into a single element, we get a ``linear''
dyadic model, that turns out to be one of the first shell models, the one
introduced by Desnianskii and Novikov in 1974~{\cite{DesNov1974JAMM}}. For
this reason we will call it DN model. Also in this case we can give a
heuristical interpretation of the model as a Littlewood-Paley decomposition:
see for example~{\cite{KatPav2005TAMS,Che2008TAMS,KisZla2005IMRN}}). The step
from the KP model to the DN one was first done by Waleffe~{\cite{Wal2006PAMS}%
}. In the same paper, he also discussed a different model, the aforementioned
Obukhov model. All three of KP, DN and Obukhov models were investigated by
Kiselev and Zlatos~{\cite{KisZla2005IMRN}}, with particular focus on the
question of regularity and blow-up.

Let us now see the DN mode in some more detail: the differential rule coupling
the intensities associated to the different shells takes the following form,
\begin{equation*}
\dot{X}_{j} = - \nu l_{j}^{2} X_{j} + c_{j - 1} X_{j - 1}^{2}
- c_{j} X_{j} X_{j + 1}, 
\end{equation*}
with $c_{j} = 2^{\alpha j}$ and $l_{j} = 2^{\gamma j}$, with $j$ taking value
in $J =\mathbb{N}$, so that $\bar{\jmath} = j - 1$, $\mathcal{O}_{j} = \{ j +
1 \}$, and $| j | \equiv j$.

This model, though physically less appealing than the KP one, has a much
simpler structure. For this reason many results, in particular regarding
uniqueness and regularity of solutions, haven been proven first for the DN
model and extended to the KP model only later.

In order to talk about existence and uniqueness of solutions, we need to state
what notion of solution are we considering for such models. A
\emph{componentwise solution} of the KP model is a family $(X_{j})_{j \in J}$
of differentiable functions such that \eqref{eq:tree_dyadic} is satisfied. If
a componentwise solution is in $L^{\infty} (\mathbb{R}_{+}, l^{2} (J))$, it is
called a \emph{Leray solution}. Analogous definitions hold for the DN model
(actually, we just have to consider $J =\mathbb{N}$, and the other conventions
written above).

\begin{theorem}
For the KP model \eqref{eq:tree_dyadic}, for any initial condition in $l^{2}$,
there exists a Leray solution.
\end{theorem}

The argument for the proof is quite standard, using Galerkin approximations,
and can be found in~{\cite{BarBiaFlaMor2013JMP,Bia2013,BiaMor2017CMP}}. A
similar result holds for the DN model, and actually, with some assumptions on
the coefficients, solutions of the DN model can be lifted to the KP model.

A natural question that can be raised at this point is the following: what
about more regular solutions? This question is in fact strongly tied to
another interesting property of dyadic models, that of anomalous dissipation.
As a matter of fact, ignoring the dissipative term, one can show that for both
KP and DN is formally preserved. However, if we approach the issue rigorously,
we see that this is only true for solutions that are regular enough. However,
it is possible to show that energy actually dissipates, hence solutions cannot
be that regular. This kind of argument is presented in the
aforementioned~{\cite{BarFlaMor2011TAMS}} for the DN model. The same is true
also for the KP model and the RCM, as it is shown
in~{\cite{BarBiaFlaMor2013JMP}} and~{\cite{BiaMor2017CMP}}.

\subsection{Uniqueness and Regularization by noise for dyadic models}


For the DN model there is uniqueness if we restrict ourselves to non-negative
solutions but we lose it if we allow for solutions that change
sign~{\cite{BarFlaMor2010CRM,BarMor2013NDEA}}. For the KP model, uniqueness
for non-negative solutions is an open problem, but counterexamples to
uniqueness can be shown for solutions that are allowed to change sign.

In both cases, counterexamples can be constructed through self-similar
solutions, that is solutions of the form $X_{j} (t) = \frac{a_{j}}{t - t_{0}}%
$, for some $t_{0} < 0$, for all $j \in J$ and $t \geqslant0$, with the
coefficient $a_{j}$ satisfying some coupling conditions. Once we have such
solutions, we can use time reversal to have solutions that blow up in finite
time, that are in particular non Leray, hence showing non-uniqueness of
componentwise solutions.

For the RCM it is hard to prove results for general solutions, because they
are quite complicated to deal with. However if we focus on constant solutions,
we not only have an existence and uniqueness result of a (finite energy)
forced solution that dissipates energy, but we can also write such solution
explicitly. In the case of the RCM, this allows us to obtain some interesting
results regarding the structure function and the geometry of the anomalous
dissipation. Existence and uniqueness of constant solutions hold for the KP
and the DN models as well~{\cite{BarBiaFlaMor2013JMP,CheFriPav2007JMP}}. In
the case of the inviscid DN model, the constant solution is particularly
interesting, because it has been proven to be a global
attractor~{\cite{CheFriPav2009D}}. A similar result holds for the viscous
model, too~{\cite{CheFri2009PDNP}}. For the KP and the RCM the existence of
such a global attractor is still a conjecture.

For constant solutions we have a uniqueness result. However this is not the
case if we consider generic solution, as mentioned above. In order to recover
some kind of uniqueness, we resort to regularization by noise techniques.

It is true that we started with a PDE, but the model that we are considering
is now made of an infinite system of coupled ODEs. So it should not be
surprising that we can obtain regularization results by adding noise. Let us
see some more details. Notice that we focus only on the inviscid case.

In order to recover uniqueness of the solution, we want to define a stochastic
perturbation of the deterministic KP model: among the several options
possible, we choose a multiplicative term (so that the perturbation ``scales''
with the solution, being neither irrelevant nor dominant) such that the total
energy is (formally) $\mathbb{P}$-a.s.~preserved
\begin{equation}
\label{eq:stoch_stratonovich}\mathrm{d} X_{j} = \left(  c_{j} X_{\bar{\jmath}%
}^{2} - X_{j} \sum_{k \in\mathcal{O}_{j}} c_{k} X_{k} \right)  \mathrm{d} t +
c_{j} X_{\overline{\jmath}} \circ\mathrm{d} W_{j} - \sum_{k \in\mathcal{O}%
_{j}} c_{k} X_{k} \circ\mathrm{d} W_{k},
\end{equation}
with $(W_{j})_{j \in J}$ a family of independent Brownian motions, together
with deterministic initial conditions $X (0) = x = (x_{j})_{j \in J} \in
l^{2}$.

For this model (which we can also write in It{\^o} formulation) we consider
solutions that are weak in the probabilistic sense. Of particular interest,
for obvious physical reasons, are energy controlled solutions, that is weak
solutions that satisfy
\[
\mathbb{P} \left(  \sum_{j \in J} X^{2}_{j} (t) \leqslant\sum_{j \in J}
x_{j}^{2} \right)  = 1 \text{for all } t \geqslant0,
\]
that is, the energy is almost surely bounded by the initial one.

\begin{theorem}
\label{thm:weak_ex_uniq_tree}There exists an energy controlled solution to
\eqref{eq:stoch_stratonovich} in $L^{\infty} (\Omega\times[0, T] ; l^{2})$ for
initial conditions $X (0) = x = (x_{j})_{j \in J} \in l^{2}$.

Moreover, there is uniqueness in law in the same class of energy controlled solutions.
\end{theorem}

Both weak existence and weak uniqueness are achieved through Girsanov theorem,
transforming the nonlinear SDEs in linear ones. However the first step is to
translate our model from the Stratonovich
formulation~\eqref{eq:stoch_stratonovich} into Ito formulation, which is
easier to manipulate:
\begin{equation}
\label{eq:stoch_ito}\mathrm{d} X_{j} = \left(  c_{j} X_{\bar{\jmath}}^{2} -
X_{j} \sum_{k \in\mathcal{O}_{j}} c_{k} X_{k} \right)  \mathrm{d} t + c_{j}
X_{\overline{\jmath}} \mathrm{d} W_{j} - \sum_{k \in\mathcal{O}_{j}} c_{k}
X_{k} \mathrm{d} W_{k} - \frac{1}{2} \left(  c_{j}^{2} + \sum_{k
\in\mathcal{O}_{j}} c_{k}^{2} \right)  X_{j} \mathrm{d} t.
\end{equation}
Moreover, since we want to use Girsanov theorem, it makes sense to
rewrite~\eqref{eq:stoch_ito} in the following form:
\[
\mathrm{d} X_{j} = c_{j} X_{\bar{\jmath}} (X_{\overline{\jmath}} \mathrm{d} t
+ c_{j} \mathrm{d} W_{j}) - \sum_{k \in\mathcal{O}_{j}} c_{k} X_{k} (X_{j}
\mathrm{d} t + \mathrm{d} W_{k}) - \frac{1}{2} \left(  c_{j}^{2} + \sum_{k
\in\mathcal{O}_{j}} c_{k}^{2} \right)  X_{j} \mathrm{d} t,
\]
where we isolated the terms $X_{\overline{\jmath}} \mathrm{d} t + c_{j}
\mathrm{d} W_{j}$, which are (for all $j \in J$) Brownian motions with respect
to a new measure $\tilde{\mathbb{P}}$ on $(\Omega, \mathcal{F}_{\infty})$.
More precisely we can state:

\begin{proposition}
Given an energy controlled solution $(\Omega, (\mathcal{F}_{t})_{t},
\mathbb{P}, W, X)$ of \eqref{eq:stoch_ito} (or equivalently
\eqref{eq:stoch_stratonovich}), we can define a measure $\tilde{\mathbb{P}}$
as follows:
\[
\left.  \frac{\mathrm{d} \tilde{\mathbb{P}}}{\mathrm{d} \mathbb{P}} \right|
_{\mathcal{F}_{t}} = \exp\left(  - \sum_{j \in J} \int_{0}^{t} X_{\bar{\jmath
}} (s) \mathrm{d} W_{j} (s) - \frac{1}{2} \int_{0}^{t} \sum_{j \in J}
X_{\bar{\jmath}}^{2} (s) \mathrm{d} s \right)  .
\]
Then the processes
\[
B_{j} (t) = W_{j} (t) + \int_{0}^{t} X_{\bar{\jmath}} (s) \mathrm{d} s
\]
are a $J$-indexed family of independent Brownian motions on $(\Omega,
(\mathcal{F}_{t})_{t}, \tilde{\mathbb{P}})$, and $(\Omega, (\mathcal{F}%
_{t})_{t}, \tilde{\mathbb{P}}, B, X)$ satisfies the linear equations
\begin{equation*}
\mathrm{d} X_{j} = c_{j} X_{\bar{\jmath}} \mathrm{d}
B_{j} (t) - \sum_{k \in\mathcal{O}_{j}} c_{k} X_{k} \mathrm{d} B_{k} (t) -
\frac{1}{2} \left(  c_{j}^{2} + \sum_{k \in\mathcal{O}_{j}} c_{k}^{2} \right)
X_{j} \mathrm{d} t. 
\end{equation*}

\end{proposition}

For this linear system we can easily prove, by Galerkin approximations, that
there exists a strong solution. The next step is to prove strong uniqueness
for the linear system.

To do so, we consider the second $\tilde{\mathbb{P}}$-moments of the $X_{j}$s:
for every solution $X$ of the (nonlinear) system~\eqref{eq:stoch_ito}, for
every $j \in J$ and $t \geqslant0$, $\tilde{\mathbb{E}} [X_{j}^{2} (t)] <
\infty$ and satisfies the differential equation:
\[
\frac{\mathrm{d}}{\mathrm{d} t} \tilde{\mathbb{E}} [X_{j}^{2} (t)] = - \left(
c_{j}^{2} + \sum_{k \in\mathcal{O}_{j}} c_{k}^{2} \right)  \tilde{\mathbb{E}}
[X_{j}^{2} (t)] + c_{j}^{2} \tilde{\mathbb{E}} [X_{\bar{\jmath}}^{2} (t)] +
\sum_{k \in\mathcal{O}_{j}} c_{k}^{2} \tilde{\mathbb{E}} [X_{k}^{2} (t)] .
\]
Now we have obtained a system of closed equations, with a very nice structure:
if we write it in matricial form, it is strongly reminiscent of the forward
equations of a Markov chain (even though it actually is not).
Thanks to this link, we can show
uniqueness for the second moments and, hence, for the solution of the linear
system. This strong uniqueness result translates into uniqueness in law for
the nonlinear system, as the two measures $\mathbb{P}$ and $\tilde{\mathbb{P}%
}$ are not equivalent on $\mathcal{F}_{\infty}$. More precise statements, as
well as detailed proofs, can be found in~{\cite{Bia2013ECP}}
and~{\cite{Bia2013}}.

A similar result holds for the DN linear dyadic model, and was obtained
earlier in~{\cite{BarFlaMor2010PAMS}}. In this case the model has the
following form:
\begin{equation*}
\mathrm{d} X_{j} = (c_{j} X_{j - 1}^{2} - c_{j + 1}
X_{j} X_{j + 1}) \mathrm{d} t + c_{j} X_{j - 1} \circ\mathrm{d} W_{j - 1} -
c_{j + 1} X_{j + 1} \circ\mathrm{d} W_{j},
\end{equation*}
with \ $(W_{j})_{j \in J =\mathbb{N}}$ a sequence of independent Brownian
motions, and the form of the noise chosen to be formally energy preserving
(almost surely). In this case, anomalous dissipation has been shown
in~{\cite{BarFlaMor2011AAP}}.

Of course one can deduce weak existence and uniqueness for DN from Theorem
\ref{thm:weak_ex_uniq_tree} for the KP model. It is interesting to notice that
the different behaviour seen in the deterministic case for non-negative and
mixed-sign solutions is now absent, even though this is not surprising,
because the noise is causing sign changes.

In the end, regularization by noise techniques had at least a partial success
in the area of fluid dynamics. Even though the techniques used for dyadic
models did not immediately translate back to Navier-Stokes equations, there are
also ideas born in the study of shell models that trickled back to
Navier-Stokes. In particular, in~{\cite{Tao2015JAMS}} some ideas from previous
works on dyadic models were used to show blow-up of an averaged version of 3D
Navier-Stokes, proving a meta-theorem: no technique that does not distinguish
the DN model from Navier-Stokes can show regularity for NSE.


\begin{thebibliography}{10}
	
	\bibitem{ArnCraWih1983SJCO}
	L.~Arnold, H.~Crauel, and V.~Wihstutz.
	\newblock Stabilization of {{Linear Systems}} by {{Noise}}.
	\newblock {\em SIAM Journal on Control and Optimization}, 21(3):451--461, 1983.
	
	\bibitem{BabMahNic1996EJMB}
	A.~Babin, A.~Mahalov, and B.~Nicolaenko.
	\newblock {Global splitting, integrability and regularity of 3D Euler and
		Navier-Stokes equations for uniformly rotating fluids}.
	\newblock {\em European Journal of Mechanics, B/Fluids}, 15(3):291--300, 1996.
	
	\bibitem{BarBesFer2014SPatA}
	D.~Barbato, H.~Bessaih, and B.~Ferrario.
	\newblock On a stochastic {{Leray}}- {$\alpha$} model of {{Euler}} equations.
	\newblock {\em Stochastic Processes and their Applications}, 124(1):199--219,
	2014.
	
	\bibitem{BarBiaFlaMor2013JMP}
	D.~Barbato, L.~A. Bianchi, F.~Flandoli, and F.~Morandin.
	\newblock A dyadic model on a tree.
	\newblock {\em Journal of Mathematical Physics}, 54(2):021507, 2013.
	
	\bibitem{BarFlaMor2010CRM}
	D.~Barbato, F.~Flandoli, and F.~Morandin.
	\newblock A theorem of uniqueness for an inviscid dyadic model.
	\newblock {\em Comptes Rendus Mathematique}, 348(9-10):525--528, 2010.
	
	\bibitem{BarFlaMor2010PAMS}
	D.~Barbato, F.~Flandoli, and F.~Morandin.
	\newblock Uniqueness for a stochastic inviscid dyadic model.
	\newblock {\em Proceedings of the American Mathematical Society},
	138(07):2607--2607, 2010.
	
	\bibitem{BarFlaMor2011AAP}
	D.~Barbato, F.~Flandoli, and F.~Morandin.
	\newblock Anomalous dissipation in a stochastic inviscid dyadic model.
	\newblock {\em The Annals of Applied Probability}, 21(6):2424--2446, 2011.
	
	\bibitem{BarFlaMor2011TAMS}
	D.~Barbato, F.~Flandoli, and F.~Morandin.
	\newblock Energy dissipation and self-similar solutions for an unforced
	inviscid dyadic model.
	\newblock {\em Transactions of the American Mathematical Society},
	363(04):1925--1925, 2011.
	
	\bibitem{BarMor2013NDEA}
	D.~Barbato and F.~Morandin.
	\newblock Positive and non-positive solutions for an inviscid dyadic model:
	Well-posedness and regularity.
	\newblock {\em Nonlinear Differential Equations and Applications NoDEA},
	20(3):1105--1123, 2013.
	
	\bibitem{BarMorRom2014AP}
	D.~Barbato, F.~Morandin, and M.~Romito.
	\newblock Global regularity for a slightly supercritical hyperdissipative
	{{Navier}}\textendash{{Stokes}} system.
	\newblock {\em Analysis \& PDE}, 7(8):2009--2027, 2014.
	
	\bibitem{Bia2013}
	L.~A. Bianchi.
	\newblock {\em Dyadic Models of Turbulence on Trees}.
	\newblock Ph.{{D}}. {{Thesis}}, {S}cuola {N}ormale {S}uperiore di {P}isa, 2013.
	
	\bibitem{Bia2013ECP}
	L.~A. Bianchi.
	\newblock Uniqueness for an inviscid stochastic dyadic model on a tree.
	\newblock {\em Electronic Communications in Probability}, 18(0):1--12, 2013.
	
	\bibitem{BiaMor2017CMP}
	L.~A. Bianchi and F.~Morandin.
	\newblock Structure {{Function}} and {{Fractal Dissipation}} for an
	{{Intermittent Inviscid Dyadic Model}}.
	\newblock {\em Communications in Mathematical Physics}, 356(1):231--260, 2017.
	
	\bibitem{BrzCapFla1992SAA}
	Z.~Brze{\'z}niak, M.~Capi{\'n}ski, and F.~Flandoli.
	\newblock Stochastic {{Navier}}-stokes equations with multiplicative noise.
	\newblock {\em Stochastic Analysis and Applications}, 10(5):523--532, 1992.
	
	\bibitem{BrzFlaMau2016ARMA}
	Z.~Brze{\'z}niak, F.~Flandoli, and M.~Maurelli.
	\newblock Existence and {{Uniqueness}} for {{Stochastic 2D Euler Flows}} with
	{{Bounded Vorticity}}.
	\newblock {\em Archive for Rational Mechanics and Analysis}, 221(1):107--142,
	2016.
	
	\bibitem{Che2008TAMS}
	A.~Cheskidov.
	\newblock Blow-up in finite time for the dyadic model of the
	{{Navier}}-{{Stokes}} equations.
	\newblock {\em Transactions of the American Mathematical Society},
	360(10):5101--5120, 2008.
	
	\bibitem{CheFri2009PDNP}
	A.~Cheskidov and S.~Friedlander.
	\newblock The vanishing viscosity limit for a dyadic model.
	\newblock {\em Phys. D}, 238(8):783--787, 2009.
	
	\bibitem{CheFriPav2007JMP}
	A.~Cheskidov, S.~Friedlander, and N.~Pavlovi{\'c}.
	\newblock Inviscid dyadic model of turbulence: {{The}} fixed point and
	{{Onsager}}'s conjecture.
	\newblock {\em Journal of Mathematical Physics}, 48(6):065503, 2007.
	
	\bibitem{CheFriPav2009D}
	A.~Cheskidov, S.~Friedlander, and N.~Pavlovi{\'c}.
	\newblock An inviscid dyadic model of turbulence: {{The}} global attractor.
	\newblock {\em Discrete and Continuous Dynamical Systems}, 26(3):781--794,
	2009.
	
	\bibitem{ConKisRyzZla2008AM}
	P.~Constantin, A.~Kiselev, L.~Ryzhik, and A.~Zlato{\v s}.
	\newblock Diffusion and {{Mixing}} in {{Fluid Flow}}.
	\newblock {\em Annals of Mathematics}, 168(2):643--674, 2008.
	
	\bibitem{CriFlaHol2019JNS}
	D.~Crisan, F.~Flandoli, and D.~D. Holm.
	\newblock Solution {{Properties}} of a {{3D Stochastic Euler Fluid Equation}}.
	\newblock {\em Journal of Nonlinear Science}, 29(3):813--870, 2019.
	
	\bibitem{Da2012}
	G.~Da~Prato.
	\newblock {\em Kolmogorov {{Equations}} for {{Stochastic PDEs}}}.
	\newblock {Birkh{\"a}user}, 2012.
	
	\bibitem{DaDeb2003JdMPeA}
	G.~Da~Prato and A.~Debussche.
	\newblock Ergodicity for the {{3D}} stochastic {{Navier}}\textendash{{Stokes}}
	equations.
	\newblock {\em Journal de Math{\'e}matiques Pures et Appliqu{\'e}es},
	82(8):877--947, 2003.
	
	\bibitem{DaFla2010JoFA}
	G.~Da~Prato and F.~Flandoli.
	\newblock Pathwise uniqueness for a class of {{SDE}} in {{Hilbert}} spaces and
	applications.
	\newblock {\em Journal of Functional Analysis}, 259(1):243--267, 2010.
	
	\bibitem{DaFlaPriRoc2013AP}
	G.~Da~Prato, F.~Flandoli, E.~Priola, and M.~R{\"o}ckner.
	\newblock Strong uniqueness for stochastic evolution equations in {{Hilbert}}
	spaces perturbed by a bounded measurable drift.
	\newblock {\em The Annals of Probability}, 41(5):3306--3344, 2013.
	
	\bibitem{DaFlaRocVer2016AP}
	G.~Da~Prato, F.~Flandoli, M.~R{\"o}ckner, and A.~Y. Veretennikov.
	\newblock Strong uniqueness for {{SDEs}} in {{Hilbert}} spaces with nonregular
	drift.
	\newblock {\em The Annals of Probability}, 44(3):1985--2023, 2016.
	
	\bibitem{DaZab2002}
	G.~Da~Prato and J.~Zabczyk.
	\newblock {\em Second {{Order Partial Differential Equations}} in {{Hilbert
				Spaces}}}.
	\newblock {Cambridge University Press}, 2002.
	
	\bibitem{DebOda2006Jee}
	A.~Debussche and C.~Odasso.
	\newblock Markov solutions for the {{3D}} stochastic
	{{Navier}}\textendash{{Stokes}} equations with state dependent noise.
	\newblock {\em Journal of Evolution Equations}, 6(2):305--324, 2006.
	
	\bibitem{DebTsu2011JdMPeA}
	A.~Debussche and Y.~Tsutsumi.
	\newblock {{1D}} quintic nonlinear {{Schr{\"o}dinger}} equation with white
	noise dispersion.
	\newblock {\em Journal de Math{\'e}matiques Pures et Appliqu{\'e}es},
	96(4):363--376, 2011.
	
	\bibitem{DelFlaVin2014CPAM}
	F.~Delarue, F.~Flandoli, and D.~Vincenzi.
	\newblock Noise {{Prevents Collapse}} of {{Vlasov}}-{{Poisson Point Charges}}.
	\newblock {\em Communications on Pure and Applied Mathematics},
	67(10):1700--1736, 2014.
	
	\bibitem{DesNov1974JAMM}
	V.~N. Desnianskii and E.~A. Novikov.
	\newblock Simulation of cascade processes in turbulent flows: {{PMM}} vol. 38,
	n{$\circeq$} 3, 1974, pp. 507\textendash{}513.
	\newblock {\em Journal of Applied Mathematics and Mechanics}, 38(3):468--475,
	1974.
	
	\bibitem{Fef2006TMPP}
	C.~L. Fefferman.
	\newblock Existence and {{Smoothness}} of the {{Navier}}-{{Stokes Equation}}.
	\newblock In J.~A. Carlson, A.~Jaffe, and A.~Wiles, editors, {\em The
		{{Millennium Prize Problems}}}, pages 55--67. {Clay Mathematics Institute
		jointly with American Mathematical Society}, {Cambridge, MA, Providence RH},
	2006.
	
	\bibitem{Fer2016EJP}
	B.~Ferrario.
	\newblock Characterization of the law for {{3D}} stochastic hyperviscous
	fluids.
	\newblock {\em Electronic Journal of Probability}, 21, 2016.
	
	\bibitem{Fla2008CIME}
	F.~Flandoli.
	\newblock An {{Introduction}} to {{3D Stochastic Fluid Dynamics}}.
	\newblock In J.~M. Morel, F.~Takens, B.~Teissier, G.~Da~Prato, and
	M.~R{\"o}ckner, editors, {\em {{SPDE}} in {{Hydrodynamic}}: {{Recent
				Progress}} and {{Prospects}}}, number 1942 in Lecture {{Notes}} in
	{{Mathematics}}, pages 51--150. {Springer Berlin Heidelberg}, {Berlin,
		Heidelberg}, 2008.
	
	\bibitem{FlaGalLuo2019}
	F.~Flandoli, L.~Galeati, and D.~Luo.
	\newblock Scaling limit of stochastic {{2D Euler}} equations with transport
	noises to the deterministic {{Navier}}-{{Stokes}} equations.
	\newblock arXiv:1905.12352, 2019.
	
	\bibitem{FlaGat1995PTRF}
	F.~Flandoli and D.~Gatarek.
	\newblock Martingale and stationary solutions for stochastic
	{{Navier}}-{{Stokes}} equations.
	\newblock {\em Probability Theory and Related Fields}, 102(3):367--391, 1995.
	
	\bibitem{FlaGubPri2010IM}
	F.~Flandoli, M.~Gubinelli, and E.~Priola.
	\newblock Well-posedness of the transport equation by~stochastic perturbation.
	\newblock {\em Inventiones Mathematicae}, 180(1):1--53, 2010.
	
	\bibitem{FlaGubPri2011SPatA}
	F.~Flandoli, M.~Gubinelli, and E.~Priola.
	\newblock Full well-posedness of point vortex dynamics corresponding to
	stochastic {{2D Euler}} equations.
	\newblock {\em Stochastic Processes and their Applications}, 121(7):1445--1463,
	2011.
	
	\bibitem{FlaLuo2018}
	F.~Flandoli and D.~Luo.
	\newblock Convergence of transport noise to {{Ornstein}}-{{Uhlenbeck}} for {{2D
			Euler}} equations under the enstrophy measure.
	\newblock arXiv:1806.09332, to appear on AoP, 2018.
	
	\bibitem{FlaLuo2019}
	F.~Flandoli and D.~Luo.
	\newblock High mode transport noise improves vorticity blow-up control in {{3D
			Navier}}-{{Stokes}} equations.
	\newblock arXiv:1910.05742, 2019.
	
	\bibitem{FlaMauNek2014JMFM}
	F.~Flandoli, M.~Maurelli, and M.~Neklyudov.
	\newblock Noise {{Prevents Infinite Stretching}} of the {{Passive Field}} in a
	{{Stochastic Vector Advection Equation}}.
	\newblock {\em Journal of Mathematical Fluid Mechanics}, 16(4):805--822, 2014.
	
	\bibitem{FlaOli2018JEE}
	F.~Flandoli and C.~Olivera.
	\newblock Well-posedness of the vector advection equations by stochastic
	perturbation.
	\newblock {\em Journal of Evolution Equations}, 18(2):277--301, 2018.
	
	\bibitem{FlaRom2002TAMS}
	F.~Flandoli and M.~Romito.
	\newblock Partial regularity for the stochastic {{Navier}}-{{Stokes}}
	equations.
	\newblock {\em Transactions of the American Mathematical Society},
	354(6):2207--2241, 2002.
	
	\bibitem{FlaRom2008PTRF}
	F.~Flandoli and M.~Romito.
	\newblock Markov selections for the {{3D}} stochastic
	{{Navier}}\textendash{{Stokes}} equations.
	\newblock {\em Probability Theory and Related Fields}, 140(3):407--458, 2008.
	
	\bibitem{Gal2019}
	L.~Galeati.
	\newblock On the convergence of stochastic transport equations to a
	deterministic parabolic one.
	\newblock arXiv:1902.06960, 2019.
	
	\bibitem{GasGes2019PTRF}
	P.~Gassiat and B.~Gess.
	\newblock Regularization by noise for stochastic
	{{Hamilton}}\textendash{{Jacobi}} equations.
	\newblock {\em Probability Theory and Related Fields}, 173(3):1063--1098, 2019.
	
	\bibitem{GesMau2018CPDE}
	B.~Gess and M.~Maurelli.
	\newblock Well-posedness by noise for scalar conservation laws.
	\newblock {\em Communications in Partial Differential Equations},
	43(12):1702--1736, 2018.
	
	\bibitem{Gyo1998SPatA}
	I.~Gy{\"o}ngy.
	\newblock Existence and uniqueness results for semilinear stochastic partial
	differential equations.
	\newblock {\em Stochastic Processes and their Applications}, 73(2):271--299,
	1998.
	
	\bibitem{Hol2015PotRSAMPaES}
	D.~D. Holm.
	\newblock Variational principles for stochastic fluid dynamics.
	\newblock {\em Proceedings of the Royal Society A: Mathematical, Physical and
		Engineering Sciences}, 471(2176):20140963, 2015.
	
	\bibitem{KatPav2005TAMS}
	N.~H. Katz and N.~Pavlovi{\'c}.
	\newblock Finite time blow-up for a dyadic model of the {{Euler}} equations
	equations.
	\newblock {\em Transactions of the American Mathematical Society},
	357(2):695--708 (electronic), 2005.
	
	\bibitem{KisZla2005IMRN}
	A.~Kiselev and A.~Zlatos.
	\newblock On {{Discrete Models}} of the {{Euler Equation}}.
	\newblock {\em International Mathematics Research Notices}, 2005:25, 2005.
	
	\bibitem{KryRoc2005PTRF}
	N.~Krylov and M.~R{\"o}ckner.
	\newblock Strong solutions of stochastic equations with singular time dependent
	drift.
	\newblock {\em Probability Theory and Related Fields}, 131(2):154--196, 2005.
	
	\bibitem{MajTimEij2001CPAM}
	A.~J. Majda, I.~Timofeyev, and E.~V. Eijnden.
	\newblock A mathematical framework for stochastic climate models.
	\newblock {\em Communications on Pure and Applied Mathematics}, 54(8):891--974,
	2001.
	
	\bibitem{MikRoz2005AP}
	R.~Mikulevicius and B.~L. Rozovskii.
	\newblock Global {{L2}}-solutions of stochastic {{Navier}}\textendash{{Stokes}}
	equations.
	\newblock {\em The Annals of Probability}, 33(1):137--176, 2005.
	
	\bibitem{MolRuzSok1985SPU}
	S.~A. Molchanov, A.~A. Ruzma{\u \i}kin, and D.~D. Sokolov.
	\newblock Kinematic dynamo in random flow.
	\newblock {\em Soviet Physics Uspekhi}, 28(4):307, 1985.
	
	\bibitem{Rom2008JSP}
	M.~Romito.
	\newblock Analysis of {{Equilibrium States}} of {{Markov Solutions}}
	to~the~{{3D Navier}}-{{Stokes Equations Driven}} by {{Additive Noise}}.
	\newblock {\em Journal of Statistical Physics}, 131(3):415--444, 2008.
	
	\bibitem{StrVar2007}
	D.~W. Stroock and S.~R.~S. Varadhan.
	\newblock {\em Multidimensional {{Diffusion Processes}}}.
	\newblock {Springer}, 2007.
	
	\bibitem{Tao2015JAMS}
	T.~Tao.
	\newblock Finite time blowup for an averaged three-dimensional
	{{Navier}}-{{Stokes}} equation.
	\newblock {\em Journal of the American Mathematical Society}, 29(3):601--674,
	2015.
	
	\bibitem{Tem1995}
	R.~Temam.
	\newblock {\em Navier-{{Stokes Equations}} and {{Nonlinear Functional
				Analysis}}: {{Second Edition}}}.
	\newblock {SIAM}, 1995.
	
	\bibitem{Ver1981MUS}
	A.~Y. Veretennikov.
	\newblock On {{Strong Solutions}} and {{Explicit Formulas}} for {{Solutions}}
	of {{Stochastic Integral Equations}}.
	\newblock {\em Mathematics of the USSR-Sbornik}, 39(3):387--403, 1981.
	
	\bibitem{Wal2006PAMS}
	F.~Waleffe.
	\newblock On some dyadic models of the {{Euler}} equations.
	\newblock {\em Proceedings of the American Mathematical Society},
	134(10):2913--2922, 2006.
	
	\bibitem{ZhuZhu2015JoDE}
	R.~Zhu and X.~Zhu.
	\newblock Three-dimensional {{Navier}}\textendash{{Stokes}} equations driven by
	space\textendash{}time white noise.
	\newblock {\em Journal of Differential Equations}, 259(9):4443--4508, 2015.
	
\end{thebibliography}
\end{document}